\input AHTOH-E.STY
\headline{\ifnum\count0=1 \vtt UDC 512.543.7+512.543.16\hss\else\hss\fi}

\centerline{\ssdbf FREE SUBGROUPS OF ONE-RELATOR RELATIVE PRESENTATIONS%
\footnote{}{\rm
This work was supported by the Russian Foundation for Basic Research,
project no. 05-01-00895.} }

\smallskip
\centerline{\ss Anton A. Klyachko}
\smallskip
{
\ssqi
\centerline{Faculty of Mechanics and Mathematics, Moscow State University}
\centerline{Moscow 119992, Leninskie gory, MSU}
\centerline{klyachko@daniil.math.msu.su}
}
\medskip
{\narrower\small\noindent
Suppose that $G$ is a nontrivial torsion-free group and $w$ is a word over
the alphabet $G\cup\{x_1^{\pm1},\dots,x_n^{\pm1}\}$.  It is proved that
for $n\ge2$ the group $\~G=\gp{G,x_1,x_2,\dots,x_n\ |\ w=1}$ always
contains a nonabelian free subgroup. For $n=1$ the question about
the existence of nonabelian free subgroups in $\~G$ is answered completely in
the unimodular case (i.e., when the exponent sum of $x_1$ in $w$ is one).
Some generalisations of these results are discussed.


\noindent{\it MSC\/\rm:} 20F05, 20E06, 20E07.

}
\medskip

\s 0. Introduction

The following theorem was stated in [Ma32]; as far as we know, the proof
appeared for the first time in~[Mo69].

\proclaim{Free subgroup theorem for one-relator groups}.
A one-relator group
$$
\gp{x_1,x_2,\dots,x_n\ |\ w=1}
$$
contains no nonabelian free subgroups if and only if it is either cyclic
or isomorphic to the Baumslag--Solitar group
$G_{1,k}=\gp{x,y\ |\ y^{-1}xy=x^k}$ for some $k\in\Z\setminus\0$.

We try to solve the same problem for relative one-relator presentations,
i.e., for groups of the form
$$
\~G=\gp{G,x_1,x_2,\dots,x_n\ |\ w=1}\:=(G*F(x_1,x_2,\dots,x_n))/\nc w.
$$
Here $G$ is a group and $w$ is an element of the free product of
$G$ and the free group $F(x_1,\dots,x_n)$. In this paper, we assume that
the group $G$ is torsion-free.

In the case $n\ge2$, the answer is as expected.

\Th 1.
If $G$ is a nontrivial torsion-free group and $n\ge2$, then the group
$\~G=\gp{G,x_1,x_2,\dots,x_n\ |\ w=1}$
contains a nonabelian free subgroup.

Note that the existence of free subgroups in $\~G$ for $n\ge3$
follows immediately from the free subgroup theorem for
one-relator groups. Thus, Theorem 1 is nontrivial only for $n=2$.

The most difficult case is $n=1$.
An important role in this situation is played by the exponent sum of the
generator in the relator. A word
$w=\prod g_it^{\epsilon_i}\in G*\gp t_\infty$ is called {\it unimodular}
if $\sum \epsilon_i=1$.

If the exponent sum of the generator in the word $w$ equals to any
number $p\ne\pm1$, then it is unknown even whether the group $G$ embeds
into $\~G=\gp{G,t\ |\ w=1}$; in other words, it is unknown when the
group $\~G$ is different from $\Z/p\Z$. There are a lot of papers on this
subject, but the answer is known only under additional strong restrictions
on the group $G$ or/and on the word $w$ (see, e.g., [B84], [KP95], [C02],
[C03], [CG95], [CG00], [EH91], [FeR98], [GR62], [IK00], [Le62], [Ly80],
[S87]).  For this reason, in this paper we study only unimodular
presentations.

The injectivity of the natural mapping $G\to\~G$ in the unimodular case
was proved in [Kl93] (see also [FeR96]). More delicate properties of the
group $\~G$ can be found in [CR01], [FoR03], [Kl06], and [Kl05].

\Th 2.
If $G$ is a torsion-free group and a word
$w\in G*\gp{t}_\infty$ is unimodular, then the group
$\~G=\gp{G,t\ |\ w=1}$
contains a nonabelian free subgroup, except in the following
two cases:
\item{\rm1)}
$w\equiv g_1tg_2$, where $g_1,g_2\in G$ (so $\~G\iso G$),
and the group $G$ contains no nonabelian free subgroups;
\item{\rm2)}
the group $G$ is cyclic and $\~G$ is isomorphic to the Baumslag--Solitar
group $G_{1,2}=\gp{g,t\ |\ g^{-1}tg=t^2}$.

In [Kl06], we suggested a generalisation of the notion of
unimodularity to the case when the word $w$ is an element of the free
product of a group $G$ and any (not necessarily cyclic) group~$T$.
A word $w\equiv g_1t_1\dots g_nt_n\in G*T$ is called {\it unimodular} if
\item{1)}
$\prod t_i$ is an element of infinite order in the group $T$;
\item{2)}
the cyclic subgroup $\gp{\prod t_i}$ is normal in $T$;
\item{3)}
the quotient group $T/\gp{\prod t_i}$ is a group with the strong
unique-product property.

Recall that a group $H$ is called a {\it UP-group, {\rm or a group with
the} unique product property}, if the product $XY$ of any two finite
nonempty subsets $X,Y\subseteq H$ contains at least one element, which
decomposes uniquely into the product of an element from $X$ and an element
from $Y$.  Some time ago, there was the conjecture that any torsion-free
group is UP (the converse is, obviously, true).  However, it turned out
that there exist counterexamples ([P88], [RS87]).

We say that a group $H$ has the {\it strong unique product property} if
the product $XY$ of any two finite
nonempty subsets $X,Y\subseteq H$ such that $|Y|\ge2$ contains at least
two uniquely decomposable elements $x_1y_1$ and $x_2y_2$ such that
$x_1,x_2\in X$,\ \ $y_1,y_2\in Y$, and $y_1\ne y_2$.

As far as we know, all known examples of UP-groups have the strong
UP-property. In particular, right orderable groups, locally indicable
groups, diffuse groups in the sense of Bowditch have the strong UP
property.

The following theorem, on the one hand, generalises (to be more precise,
complements) Theorem 2, and on the other hand, is a key step in the proof
of Theorem 1.

\Th 3.
If a group $G$ is torsion-free, a group $T$ is noncyclic, and a word
$w\in G*T$ is unimodular, then the group
$\~G=\gp{G,T\ |\ w=1}\:=(G*T)/\nc w$
contains no nonabelian free subgroups if and only if
$G$ is cyclic, $T$ contains no nonabelian free
subgroups, and $w$ is conjugate in $G*T$ to a word of the form $gt$, where
$t\in T$ and $g$ is a generator of~$G$.

\proclaim{The notation} \rm which we use is mainly standard. \rm
If $x$ and $y$ are elements of a group, and $X$ is a subset
of this group, then $x^y$ means $y^{-1}xy$, commutator $[x,y]$ is
understood as $x^{-1}y^{-1}xy$, and the symbols $\gp{X}$ and $\nc{X}$
denote, respectively, the subgroup generated by the set $X$ and the
normal subgroup generated by the set $X$.
The symbol $|X|$ denotes the cardinality of $X$.


\s 1. Proof of Theorem 1

Suppose that the word $w$ has the form
$w\equiv g_1x_{j_1}^{\epsilon_1}g_2x_{j_2}^{\epsilon_2}\dots
g_px_{j_p}^{\epsilon_p}$ and the word $w'\in F(x_1,\dots, x_n)$ is
obtained from $w$ by erasing coefficients:
$w'=x_{j_1}^{\epsilon_1}x_{j_2}^{\epsilon_2}\dots x_{j_p}^{\epsilon_p}$.

\smallskip
\noindent{\bf Case 1}:
$w'$ is a proper power in the free group $F(x_1,\dots, x_n)$.
In this case, $\~G$ contains a nonabelian free subgroup,
because, according to the free subgroup theorem for one-relator groups,
a nonabelian free subgroup exists in the one-relator group
$T_1=\gp{x_1,\dots,x_n\ |\ w'=1}$, which is a homomorphic image of $G$.

\smallskip
\noindent{\bf Case 2}:
$w'$ is not a proper power. Consider the groups
$$
T=\gp{x_1,\dots,x_n\ |\ [x_1,w']=\dots=[x_n,w']=1}
\quad\hbox{and}\quad
T_1=\gp{x_1,\dots,x_n\ |\ w'=1}=T/\nc{w'}.
$$
The group $T$ is the free central extension of the one-relator group
$T_1$. It is well known that, if $w'$ is not a proper power in the free
group $F(x_1,\dots,x_n)$, then the group $T_1$ is locally indicable ([B84])
and, therefore, has the strong unique product property. The element
$w'$ has infinite order in the group $T$ (see [LS77]). Thus, the word $w$
considered as an element of the free product $G*T$ is unimodular.
The group $T$ is noncyclic, because its commutator quotient
is a free abelian group of rank $n$ and $n\ge2$.
It remains to note that the group
$\gp{G,T\ |\ w=1}$ is a homomorphic image of $\~G$ and, thus, the
assertion of Theorem 1 follows immediately from Theorem~3.


\s 2. Proof of Theorem 2

Following [FoR03], we say that an element $v\in G*\gp t_\infty$ with
cyclically reduced form
$v\equiv
g_1t^{\epsilon_1}\dots g_nt^{\epsilon_n}$, where
$\epsilon_i\in\{\pm1\}$ and $g_i\in G$, has {\it complexity} 0 if all
exponents $\epsilon_i$ are equal (i.e., the word is either
positive or negative); we say that the complexity of the word $v$ is 1 if
the cyclic sequence $(\epsilon_1,\dots,\epsilon_n)$ contains both
positive and negative exponents, but either two successive exponents are
never both positive, or two successive exponents are never both negative.
In all other cases, we say that the complexity of $v$ is higher than 1.
The complete definition can be found in [FoR03].

\proclaim{Minimal complexity theorem} {\rm [FoR03]}.
If a group $G$ is torsion-free, a cyclically reduced word $w\in G*\gp
t_\infty$ is unimodular and the complexity of a word $v\in G*\gp t_\infty$
is lower than that of $w$, then $v\ne1$ in the group
$\~G=\gp{G,t\ |\ w=1}$.

This theorem implies immediately that, in the case when the complexity
of $w$ is higher than 1, $\~G$ contains the free square of the group $G$
(because $\gp{G,G^t}=G*G^t$) and, hence, a nonabelian free subgroup.

It remains to consider a presentation with the relator of complexity 1:
$$
\~G=\gp{G, t\ \Biggm|\ ct\prod_{i=0}^m(b_ia_i^t)=1},
\quad\hbox{where } m\ge0,\ a_i,b_i\in G\setminus\1,\ c\in G.
\eqno{(1)}
$$
In this case, we use the following lemma.

\Lemma {\rm [Kl05, Lemma 22]}.
If a group $G$ is torsion-free and $\~G$ has presentation
$(1)$, then
there exists a $d\in\{2,3\}$ such that
$$
u\equiv\prod_{i=1}^s y_ix_i^{t^d}\ne1\ {\rm in}\ \~G
$$
for any positive integer $s$ and any $x_i, y_i\in G$
for which
$\Bigl|\{i\ |\ x_i\in\gp{a_m}\}\Bigr|+
\Bigl|\{i\ |\ y_i\in\gp{b_0}\}\Bigr|\le2$ and $u\ne1$ in $G*\gp{t}_\infty$.


This lemma implies that the elements $g_1h_1^{t^d}$ and $h_2^{t^d}g_2$ of
the group $\~G$ generate a free subgroup of rank 2 for any
$g_1,g_2,h_1,h_2\in G$ such that
$h_1,h_2,h_1h_2\notin\gp{a_m}$ and $g_2,g_1,g_2g_1\notin\gp{b_0}$.
Therefore, the absence
of nonabelian free subgroups in $\~G$ implies that the group $G$ is
virtually cyclic ($|G:\gp{a_m}|\le2$ or $|G:\gp{b_0}|\le2$) and, hence,
$G$ is cyclic (because $G$ is torsion-free).


If the group $G$ is cyclic, then $\~G$ is a one-relator group.
The free subgroup theorem for one-relator groups
implies that either $\~G$ contains a nonabelian free subgroup,
$\~G$ is cyclic (which is imposible for $m\ge0$), or $\~G$ is isomorphic
to $G_{1,k}$. In the last case, the unimodularity of the
relation $w$ implies that the number $k$ must be 2. It is easy to verify
this by considering the commutator quotient. Theorem 2 is proven.


\s 3. Proof of Theorem 3

Suppose that the word $w$ has the form $w\equiv g_1t_1\dots g_nt_n$.

\smallskip
\noindent
{\bf Case 1}: $n=1$.
Clearly, in this case $|T/\gp{t_1}|=\infty$,
$$
\~G\iso\cases{
  G\zvezda_{\kern 5pt\raise10pt\hbox{${}^{g_1=t_1^{-1}}$}} T
                  \quad&\hbox{if $g_1\ne1$ in $G$},\cr
  G*(T/\gp{t_1}) \quad&\hbox{if $g_1=1$ in $G$,}
   }
$$
and the assertion of theorem is obtained from the following well-known
simple facts:
{\sl
\item{\bf 1.}
An amalgamated product contains a nonabelian
free subgroup if the amalgamated subgroup is proper in
each factor and its index is larger than 2 in one of the factors.
\item{\bf 2.}
Let $\gp a$ be a cyclic normal subgroup of a group $A$. Then $A$
contains a nonabelian free subgroup if and only if such a
subgroup exists in the quotient group $A/\gp a$.

}

\smallskip
\noindent
{\bf Case 2}:
$n>1$ and the group $\gp{t_1,\dots,t_n}$ is cyclic (and,
therefore, is generated by the element $t=\prod t_i$ by virtue of
unimodularity). In this case, $\~G$ is the
amalgamted product:
$$
\~G\iso\gp{G,t\ |\
w=1}\zvezda_{\gp t} T.
$$
This free product always contains a nonabelian free
subgroup, because the amalgamated subgroup has infinite index in
both factors. The infinity of the index of $\gp t$ in the first
factor follows from the minimal complexity theorem, which, in particular,
claims that $G\cap\gp t =\1$ for $n>1$. The equality $|T:\gp t|=\infty$
holds because $w$ is unimodular and $T$ is noncyclic.

\smallskip
\noindent
{\bf Case 3}:
The group $\gp{t_1,\dots,t_n}$ is noncyclic.
In this case, without loss of generality we can assume that
$T=\gp{t_1,\dots,t_n}$.

Put $t=\prod t_i$. Let us decompose $T$ into the union of cosets:
$$
T=\coprod_{x\in T/\gp{t}} c_x\gp{t}, \quad\hbox{where }c_1=1,
$$
and rewrite the word $w$ in the form
$$
w\equiv t\prod_i {g_i}^{c_{x_i}t^{k_i}}=1.
\eqno{(2)}
$$
Let $X_1=\{x_i\}$ be the set of all $x\in T/\gp{t}$ occurring in
the reduced expression (2). Note that $|X_1|>1$, because the group
$T=\gp{t_1,\dots,t_n}$ is noncyclic.
In [Kl06], it was shown that in the group~$\~G$ we have the decomposition
$$
H_1=\gp{\{G^{c_y}\ |\ y\in X_1\}}=\zvezda_{y\in X_1}G^{c_y}.
$$
This implies immediately that $\~G$ contains the free square
of the group $G$ and, hence, a nonabelian free subgroup.


\s{\rm REFERENCES}

\medskip

\item{[B84]}
Brodskii S. D.
Equations over groups and one-relator groups
{// Sib. Mat. Zh.} 1984. {T.25}. no.2. P.84--103.

\item{[KP95]}
Klyachko Ant. A., Prishchepov M. I.
The descent method for equations over  groups
{// Moscow Univ. Math. Bull.} 1995, {V.50}  P. 56--58.

\item{[Kl05]}
Klyachko Ant. A.
The Kervaire--Laudenbach conjecture and presentations of simple groups
{// Algebra and Logic.} 2005. {T. 44}. {no.4}. P. 219--242.

\item{[Kl06]}
Klyachko Ant. A.
How to generalize known results on equations over groups
{// Mat. Zametki.} (to appear).
See also arXiv:math.GR/0406382.

\item{[LS77]}
Lyndon R. C., Schupp P. E.
{Combinatorial group theory},
Springer-Verlag, Berlin/Heidelberg/New~York, 1977.

\item{[Mo69]}
Moldavanskii D. I.
On a theorem of Magnus (in Russian)
{// Uch. zap. Ivanovsk. gos. ped. inst.} 1969. {T.44}. S.26--28.

\item{[C02]}
Clifford A.
A class of exponent sum two equations over groups
{// Glasgow Math. J.} 2002. {V.44}. P.201--207.

\item{[C03]}
Clifford A.
Nonamenable type K equations over groups
{// Glasgow Math. J.} 2003. {V.45}. P.389--400.

\item{[CG95]}
Clifford A., Goldstein R. Z.
Tesselations of $S^2$ and equations over torsion-free groups
{// Proc. Edinburgh Math. Soc.} 1995. {V.38}. P.485--493.

\item{[CG00]}
Clifford A., Goldstein R. Z.
Equations with torsion-free coefficients
{// Proc. Edinburgh Math. Soc.} 2000. {V.43}. P.295--307.

\item{[CR01]}
Cohen M. M., Rourke C.
The surjectivity problem for one-generator, one-relator extensions of
torsion-free groups
{// Geometry \& Topology}. 2001. {V.5}. P.127--142.

\item{[EH91]}
Edjvet M., Howie J.
The solution of length four equations over groups
{// Trans. Amer. Math. Soc.} 1991. {V.326}. P.345--369.

\item{[FeR96]}
Fenn R., Rourke C.
Klyachko's methods and the solution of equations over torsion-free groups
{// L'Enseignment Math\'ematique.} 1996. {T.42}. P.49--74.

\item{[FeR98]}
Fenn R., Rourke C.
Characterisation of a class of equations with solution over torsion-free
groups,
from {``The Epstein Birthday Schrift"}
{(eds. I. Rivin, C. Rourke and C. Series)},
{Geometry and Topology Monographs.} 1998. {V.1}. P.159-166.

\item{[FoR03]}
Forester M., Rourke C.
Diagrams and the second homotopy group
{// arXiv:math.AT/0306088}.

\item{[GR62]}
Gerstenhaber M., Rothaus O. S.
The solution of sets of equations in groups
{//  Proc. Nat. Acad. Sci. USA}. 1962. {V.48} P.1531--1533.

\item{[IK00]}
Ivanov S. V., Klyachko Ant. A.
Solving equations of length at most six over torsion-free groups
{// J. Group Theory}. 2000. {V.3}. P.329--337.

\item{[Kl93]}
Klyachko Ant. A.
A funny property of a sphere and equations over groups
{// Comm. Algebra}. 1993. {V.21}. P.2555--2575.

\item{[Le62]}
Levin F.
Solutions of equations over groups
{// Bull. Amer. Math. Soc.} 1962. {V.68}. P.603--604.

\item{[Ly80]}
Lyndon R. C.
Equations in groups
{// Bol. Soc. Bras. Math.} 1980. {V.11}. \number1. P.79--102.

\item{[Ma32]}
Magnus W.
Das Identit\"atsproblem f\"ur Gruppen mit einer definierenden Relation
{// Math. Ann.} 1932. {V.106}. P.295--307.

\item{[P88]}
Promyslow S. D.
A simple example of a torsion free nonunique product group
{// Bull. London Math. Soc.} 1988. {V.20}. P.302--304.

\item{[RS87]}
Rips E., Segev Y.
Torsion free groups without unique product property
{// J. Algebra} 1987. {V.108}. P.116--126.

\item{[S87]}
Stallings J. R.
A graph-theoretic lemma and group embeddings
// Combinatorial group theory and topology
(eds. S. M. Gersten, J. R. Stallings).
{Annals of Mathematical Studies}. 1987. {V.111}. P.145--155.

\end